\documentclass[12pt]{article}
\usepackage[centertags]{amsmath}
\usepackage{amsfonts}
\usepackage{amssymb}
\usepackage{amsthm}
\input{amssym.def}

\numberwithin{equation}{section}

\newtheorem{Definition}{Definition}[section]
\newtheorem{definition}[Definition]{Definition}
\newtheorem{theorem}[Definition]{Theorem}
\newtheorem{lemma}[Definition]{Lemma}

\newtheorem{corollary}[Definition]{Corollary}

\begin{document}

\title{\Large \bf Rank preservers of matrices over additively idempotent and multiplicatively cancellative semirings}
\author{\bf A. K. Bhuniya and Sushobhan Maity}
\date{}

\maketitle
\begin{center}
Department of Mathematics, Visva-Bharati, Santiniketan-731235, India. \\
anjankbhuniya@gmail.com, susbhnmaity@gmail.com
\end{center}

\begin{abstract}{\footnotesize }
Here we characterize the linear operators that preserve rank of matrices over additively idempotent and multiplicatively cancellative semirings. The main results in this article generalize the corresponding results on the two element Boolean algebra \cite{BP} and on the max algebra \cite{Bapat}; and holds on max-plus algebra and some other tropical semirings.
\end{abstract}

\section{Introduction}
\pagenumbering{arabic}

There has been a great deal of interest in recent years in the study of linear operators that preserve the rank (there are different notions of rank) of matrices over distributive lattices \cite{BP3,SV,SP}, Boolean algebra \cite{BP,SHJK,Kang}, max algebra \cite{Bapat,SK,Song}, tropical semirings, semirings \cite{BP2,KSB, Tan} etc. Major attraction in this direction has been shown to the matrices over either on the two element Boolean algebra $\mathbb{B}_2=\{0,1\}$ or on the max algebra $\mathbb{R}_{max}=(\mathbb{R}^+ \cup \{0\},max,\cdot)$. Interestingly, all of $\mathbb{B}_2$, $\mathbb{R}_{max}$, max-plus algebra over $\mathbb{R} \cup \{-\infty \}$, tropical semirings $(\mathbb{Z}, max,+)$, $(\mathbb{N} \cup \{0\},max,\cdot)$, etc are multiplicatively cancellative semirings with an idempotent additive reduct.

In \cite{BP}, Beasley and Pullman studied the factor rank preserving linear operators on the linear spaces of $m \times n$ matrices over the two element Boolean algebra $\mathbb{B}_2$. They obtained many results analogous and near analogous to some results on linear operators on field-valued matrices. Bapat, Pati and Song further generalized these results to the linear operators on the semimodules over max algebra $\mathbb{R}_{max}$ consisting of the nonnegative real numbers equipped with two binary operations maximization and multiplication. Since the class of all additively idempotent and multiplicatively cancellative semirings contains $\mathbb{B}_2$ and $\mathbb{R}_{max}$ as well as many other interesting semirings including the max-plus algebra $\mathbb{R}^{max}$ and tropical semiring $(\mathbb{Z}, max, +)$, etc., so it is worth to study the extent up to which the results obtained in \cite{BP} and \cite{Bapat} can be carried over an additively idempotent and multiplicatively cancellative semiring.

Here we show that a linear operator $T$ on the semimodule of all $m \times n$ matrices over such a semiring is a rank preserver if and only if it preserves the rank of all rank-1 and rank-2 matrices.

\section{Preliminaries and basic results}
A \emph{semiring} \index{semifield} $(S,+, \cdot)$ is an algebra with two binary operations $'+'$ and $'\cdot '$ such that
\begin{enumerate}
\item[(i)] $(S, +)$ is a commutative monoid with zero element 0;
\item[(ii)] $(S, \cdot)$ is a commutative monoid with unit element 1;
\item[(iii)] the following distributive laws hold:
\begin{align*}
 x(y+z) = xy+xz \; \textrm{and} \; (x+y)z = xz+yz.
\end{align*}
\item[(iv)] the zero $0$ is absorbing, thai is $a \cdot 0=0 \cdot a=0$ for all $a \in S$.
\end{enumerate}

A nonzero element $a \in S$ is called a zero divisor if $ab=0$ for some nonzero $b \in S$. An element $u$ of $S$ is a unit if there exists an element $v \in S$ such that $uv=1$. The element $v$ is called the inverse of $u$ in $S$. We will denote the set of all units by $U(S)$. An element $a$ of $S$ is multiplicatively cancellable if  $ba=ca$ only when $b=c$. Clearly every unit of $S$ is multiplicatively cancellable and no multiplicatively cancellable element of $S$ is a zero divisor. If every nonzero element of $S$ is multiplicatively cancellable then we say that the semiring $S$ is multiplicatively cancellable.

For more on semirings, we refer to \cite{golan}.
\begin{definition}
A semiring $S$ is said to be additively unit irreducible if $a+b \in U(S)$ implies that either $a\in U(S)$ or $b \in U(S)$.
\end{definition}

Throughout this article, unless otherwise stated, $S$ always mean an additively unit irreducible semiring which is multiplicatively cancellable and additively idempotent, that is $a+a=a$ for all $a \in S$. Each of the semirings $\mathbb{B}_2$, $\mathbb{R}_{max}$, $\mathbb{R}^{max}$, $(\mathbb{Z}, max, +)$ are of this type.

Following observation is due to Mora, Wasanawichit, Kemprasit {\cite{MWK}}.
\begin{lemma}{\cite{MWK}}
Let $(S,+, \cdot)$ be an additively idempotent semiring with zero 0 and identity 1. Then for every $a,b \in S$, $a+b=0$ $\Rightarrow $ $a=0$ and $b=0$.
\end{lemma}

Thus every semiring $S$, we are considering here, is zerosumfree. This result is a direct consequence of additive idempotency. As a consequence of multiplicatively cancellativeness, it follows that there is no zero divisor in $S$. On $S$ a partial order $'\leq '$, is defined by: for $a,b \in S$,
\begin{center}
$a \leq b$  if  $a+b=b$.
\end{center}

Denote the set of all $m \times n$ matrices with entries from $S$ by $M_{m \times n}(S)$ and $\Delta =\{(i,j):1 \leq i \leq m; 1 \leq j  \leq n\}$. The $(i,j)$th entry of an $m\times n$ matrix $A$ over $S$ is denoted by $A_{ij}$ or $A(i,j)$ for every $(i,j) \in \Delta$. If $A,B \in M_{m \times n}(S)$ are two $m \times n$ matrices over $S$, then their product, sum, transpose and muiltiplication by a scalar are defined in the usual way. A square matrix $M$ is called a monomial if it has exactly one nonzero element in each row and column.

\begin{lemma}{\cite{SV}}[Theorem 1]
Let $S$ be a zerosumfree commutative semiring with no zero divisors. Then $A \in M_n(S)$ is invertible if and only if it is a monomial matrix all of whose nonzero entries are units.
\end{lemma}

The natural partial order on $S$ induces a partial order $'\leq '$ on $M_{m\times n}(S)$ given by: for $A=[A_{ij}]$ and $B=[B_{ij}]$ in $M_{m\times n}(S)$,
\begin{center}
$A \leq B$ if $A_{ij} \leq B_{ij}$ for all $(i,j) \in \Delta$.
\end{center}
Thus $A \leq B$ if $A_{ij}+B_{ij}=B_{ij}$ for all $(i,j) \in \Delta$. If $A \leq B$, then we say that $B$ dominates $A$.

A nonempty subset $V$ of $M_{m \times n}(S)$ is said to be a semimodule if $V$ is closed under addition and scalar multiplication. If $V$ and $W$ are semimodules over $S$ with $V \subseteq W$, then $V$ is called a subsemimodule of $W$. The semimodule of all $1 \times n$ matrices over $S$ is of special interest, which we denote by $S^n$.

Let $G$ be a subset of $M_{m \times n}(S)$. Then $span(G)$ is defined by
\begin{center}
$span(G)=\{x: x=\sum ^k_{i=1} \alpha _ix_i |k \in \mathbb{N}, \alpha _i \in S, x_i \in G\}$.
\end{center}
If there is a finite subset $G$ such that $V=span(G)$, then $V$ is called finitely generated.

Let $V$ be a semimodule over a semiring $S$. A set $D$ of vectors in $V$ is called linearly dependent if there exists $x \in D$ such that $x \in span(D - \{x\})$; otherwise it is called linearly independent. Thus a linearly independent set cannot contain the zero vector. A subset $B$ of $V$ is called a basis of $V$ if $span(B)=V$ and $B$ is linearly independent. Every finitely generated semimodule has a finite basis.

Denote by $E_{ij}$ the $m \times n$ matrix over $S$ such that the $(i,j)$ th entry is 1 and all other entries are 0. Then $\mathbb{E}=\{E_{ij}:(i,j) \in \Delta \}$ is a basis of the semimodule $M_{m \times n}(S)$, which we call the standard basis of the semimodule $M_{m \times n}(S)$.

\begin{theorem}       \label{equality of basis}
Let $V$ be a semimodule over a semiring $S$. Let $\mathcal{B}_1,\mathcal{B}_2$ be two bases of $V$. Then for $x \in \mathcal{B}_1$, there exists a unique $y_x \in \mathcal{B}_2$ such that $y_x=\alpha x$ for some unit $\alpha \in S$ and for $y \in \mathcal{B}_2$ there exists a unique $x_y \in \mathcal{B}_1$ such that $x_y=\beta y$ for some unit $\beta \in S$. In particular $| \mathcal{B}_1 |= |\mathcal{B}_2 |$.
\begin{proof}
Let $x \in \mathcal{B}_1$, then there exist $c_1,c_2, \cdots ,c_n \in S$ and vectors $y_1,y_2, \cdots ,y_n \in \mathcal{B}_2$ such that
\begin{center}
$x=\sum^n_{i=1} c_iy_i$   $\;\;\; \cdots \cdots (1)$
\end{center}
Since each $y_i$ is a linear combination of a finite numbers of elements of $\mathcal{B}_1$, there exist $x_1,x_2, \cdots ,x_m \in \mathcal{B}_1$ such that
\begin{center}
$y_i=\sum^m_{j=1}a^j_ix_j$,  $\;\;\; \cdots \cdots (2)$
\end{center}

From (1) and (2) we have,
\begin{align*}
x=&\sum^n_{i=1} c_i\sum^m_{j=1}a^j_ix_j
\\& =\sum^n_{i=1}c_ia^1_ix_1 + \sum^n_{i=1}c_ia^2_ix_2 + \cdots + \sum^n_{i=1}c_ia^m_ix_m.
\end{align*}

Since $\mathcal{B}_1$ is a basis of $V$ and $x\in \mathcal{B}_1$, it follows that $x=x_{j_0}$ for some $j_0$. Thus $\sum^n_{i=1}c_ia^{j_0}_i=1$ and $c_ia^j_i=0$ for all $j \neq j_0$ and for all $i$. From this we see that, there exists $i_0$ such that $c_{i_0}a^{j_0}_{i_0}$ is a unit and hence $a^{j_0}_{i_0}$ is a unit. So $c_{i_0} \neq 0$ and $a^{j}_{i_0}=0$ for all $j \neq j_0$. Thus $a^j_{i_0}x_j \leq a^{j_0}_{i_0}x_{j_0}$, for all $j\neq j_0$. From (2), $y_{i_0}=\sum^m_{j=1}a^j_{i_0}x_j=a^{j_0}_{i_0}x_{j_0}=a^{j_0}_{i_0}x$.

Thus for $x \in \mathcal{B}_1$, there exists a $y_x \in \mathcal{B}_2$ i.e. $y_{i_0}$ and a unit $\alpha=a^{j_0}_{i_0}$ such that $y_x=\alpha x$. If for $y_1,y _2 \in \mathcal{B}_2$, there are $\alpha_1, \alpha_2 \in S$ such that $y_1=\alpha_1x$ and $y_2=\alpha_2x$, then there exists a unit $\gamma \in S$ such that $y_1=\gamma y_2$, which is a contradiction.
Similarly it can be shown that for each $y\in \mathcal{B}_2$ there exists a unique $x_y \in \mathcal{B}_1$ such that $x_y=\beta y$, for some unit $\beta \in S$.

The above discussion shows that the function $f: \mathcal{B}_1 \longrightarrow \mathcal{B}_2$ defined by $f(x)=y_x$ is a bijection. Thus the proof is complete.
\end{proof}
\end{theorem}

From the above result it follows that every basis of a finitely generated semimodule over a semiring contains the same number of vectors. This number which is finite is called the dimension of $V$ and is denoted by dim$(V)$.

\begin{corollary}
Every basis $\beta =\{\beta _1, \beta _2, \cdots ,\beta _n \}$ of $S^n$ is of the form $\beta _i=\left(
                                                                                          \begin{array}{c}
                                                                                            0 \\
                                                                                            \vdots \\
                                                                                            a_i \\
                                                                                            \vdots \\
                                                                                            0 \\
                                                                                          \end{array}
                                                                                        \right)$, where $a_i \in U(S)$.
\end{corollary}

If $V,W$ are semimodules over a semiring $S$, a mapping $T:V\rightarrow W$ is called a linear transformation if $T$  has the following properties: for every $\alpha , \beta \in S$ and $x,y \in V$,
\begin{center}
$T(\alpha x+ \beta y)=\alpha T(x) + \beta T(y)$
\end{center}
Then it follows that $T(0)=0$. If $V=W$, then $T$ is called a linear operator.

If $T$ is linear , then $\alpha \leq \beta$ implies that $T(\alpha) \leq T(\beta)$. Let $V$ and $W$ are two semimodules over $S$. A linear transformation $T:V \rightarrow W$ is called injective if $T(x)=T(y)$ implies $x=y$ for all $x,y \in V$ and is called surjective if $T(V)=W$. It is called invertible if it is injective and surjective. If $\mathfrak{B}$ is a basis of $V$, then the image of $V$ in $W$, $T(V)$ is generated by the image $T(\mathfrak{B})$. Thus we have:
\begin{lemma}
Let $V$ be a finitely generated semimodule and $T$ be a linear operator on $V$. Then, dim$(T(U))\leq $dim$(U)$ for every subsemimodule $U$ of $V$.
\end{lemma}

If a linear transformation $T:V \rightarrow W$ is such that dim$(T(U))$=dim$(U)$ for every subsemimodule $U$ of $V$, then $T$ is said to preserve dimension.
\begin{lemma}      \label{injectivity}
Let $V$ and $W$ be finitely generated semimodules. If $T:V\rightarrow W$ is injective, then $T$ preserves dimension and $T$ maps every basis of $V$ onto a basis of $T(V)$.
\begin{proof}
Let $\{v_1,v_2, \cdots ,v_n \}$ be a basis of $V$. Then $\{T(v_1),T(v_2), \cdots ,T(v_n)\}$ spans $T(V)$. Let $T(v_j)=\sum^n_{i=1,i \neq j}a_iT(v_i)$, then $T$ being injective, $v_j=\sum^n_{i=1,i \neq j}a_iv_i$; which contradicts that $\{v_1,v_2, \cdots ,v_n \}$ is linearly independent. Thus $\{T(v_1),T(v_2), \cdots ,T(v_n)\}$  is a basis of $T(V)$ and the result follows.
\end{proof}

Hence we will deal with finitely generated semimodules over $S$ only.
\end{lemma}
\begin{lemma}    \label{preserv every subspace}
If $T:V \rightarrow W$ is a surjective linear transformation, then $T$ is invertible if and only if $T$ preserves dimension.
\begin{proof}
If $T$ is invertible, then $T$ is injective and so preserves dimension, by Lemma \ref{injectivity}.

Conversely, assume that $T(x)=T(y)$. Let $U$ be the subsemimodule generated by $x,y$. Since $T(x)$ generates $T(U)$, so dim$(T(U))=1$ and hence dim$(U)=1$. Then $x=\alpha A$ and $y=\beta A$ for some $\alpha , \beta \in S$. Now $\beta T(x)=\alpha T(y)$ implies that $\alpha=\beta$, by the multiplicative cancellative property of $S$. Thus $T$ is injective and hence is invertible.
\end{proof}
\end{lemma}
\begin{corollary}
Let $T$ be a linear operator on $V$. Then the following conditions are equivalent:
\begin{enumerate}      \label{invertible of t}
\item[(i)] $T$ is invertible;
\item[(ii)] $T$ preserves dimension;
\item[(iii)] $T$ permutes every basis of $V$, with some unit scalar multiplication;
\item[(iv)] T permutes the standard basis with some unit scalar multiplication;
\end{enumerate}
\begin{proof}
Equivalence of $(i)$ and $(ii)$ follows from Lemma \ref{preserv every subspace}. To show that $(i)$ and $(iii)$ are equivalent, note that if $T$ is invertible and $\{v_1,v_2, \cdots ,v_n\}$ is a basis of $V$, then $\{T(v_1),T(v_2), \cdots ,T(v_n)\}$ is a basis of $T(V)=V$. Thus from Theorem \ref{equality of basis}, $(i)$ implies $(iii)$. On the otherhand, let $T$ be a linear operator satisfying $(iii)$. Then for all $v \in V$, $T(v)=Mv$, where $M$ is a monomial whose nonzero entries are units, consequently $T$ is invertible. From above discussion it follows that $(i)$ and $(iv)$ are equivalent.
\end{proof}
\end{corollary}

It follows directly from Corollary \ref{invertible of t} that,
\begin{corollary}
A linear operator on $M_{m \times n}(S)$ is invertible if and only if $T$ permutes every basis of $M_{m \times n}(S)$ if and only if $T$ preserves dimension.
\end{corollary}

\section{Rank-1 preserving linear operators}     \label{rank-1 preserver}
If $A$ is an $m\times n$ matrix with entries from $S$, then the rank of $A$, denoted by $r(A)$, is the least positive integer k for which there exist $m \times k$ and $k \times n$ matrices $B$ and $C$ over $S$ such that $A=BC$. The rank of the zero matrix is $0$.

It is easy to check that $r(A)$ is the least positive integer k such that $A$ can be expressed as a sum of k matrices of rank-1.

A linear operator $T$ on $M_{m \times n}(S)$ is called
\begin{enumerate}
\item[(i)] a $(U,V)$ operator if there exist invertible matrices $U$ of order $m$ and $V$ of order $n$, such that $T(A)=UAV$ for all $A \in M_{m \times n}(S)$ or m=n and $T(A)=UA^tV$ for all $A \in M_{m \times n}(S)$.
\item[(ii)] rank preserver if $r(T(A))=r(A)$ for all $A \in M_{m \times n}(S)$.
\item[(iii)] rank-1 preserver if $r(A)=1$ implies $r(T(A))=1$ for all $A \in M_{m \times n}(S)$.
\end{enumerate}
Note that every $(U,V)$ operator is a rank preserver.

We call a subsemimodule of $M_{m \times n}(S)$ whose nonzero members have rank 1 as a rank-1 subsemimodule.
\begin{lemma}   \label{rank reuction }
If $T$ is a linear operator on $M_{m \times n}(S)$ that preserves the dimension of all rank-1 subsemimodules, then the restriction of $T$ to the set of all rank-1 matrices is injective or $T$ reduces the rank of some rank-2 matricx to 1.
\begin{proof}
Let $\mathcal{M}^1=\{A\in M_{m \times n}(S):r(A)=1\}$. For each $B \in \mathcal{M}^1$, define $\mathcal{W}_{B}=span \{X \in \mathcal{M}^1: T(X)=T(B)\}$. Note that $B \in \mathcal{W}_{B}$ and dim$(T(\mathcal{W}_{B}))=1$. Then we have two cases:

\textbf{Case 1}: If $\mathcal{W}_{B}$ is a rank-1 subsemimodule, then from hypothesis, dim$(\mathcal{W}_{B})=$ dim$(T(\mathcal{W}_{B}))=1$. So $\mathcal{W}_{B}=\langle B \rangle$. Hence $T$ is injective.

\textbf{Case 2}: There exists $B \in \mathcal{M}^1$ such that dim$(\mathcal{W}_{B}) >1$. Then there are $X,Y \in \{X \in \mathcal{M}^1: T(X)=T(B)\}$ such that $r(X+Y)=2$, but $r(T(X+Y))=r(T(X)+T(Y))=r(T(B))=1$. Hence $T$ reduces the rank of some rank-2 matrix to 1.
\end{proof}
\end{lemma}
\begin{corollary}      \label{rank-1 rank-2 invertible}
If $T$ is a linear operator on $M_{m \times n}(S)$ that
\begin{enumerate}
\item[(i)] preserves the rank of all rank-1 and rank-2 matrices and
\item[(ii)] preserves the dimension of all rank-1 subsemimodules,
\end{enumerate}
then
\begin{enumerate}
\item[(a)] $T$ is invertible and
\item[(b)] $T^{-1}$ satisfies $(i)$ and $(ii)$.
\end{enumerate}
\begin{proof}
Let $\mathbb{E}$ be the standard basis of $M_{m \times n}(S)$. To prove that T is invertible it is sufficient to show that $T$ permutes $\mathbb{E}$ with some unit scalar multiplication, by Corollary \ref{invertible of t}. Lemma \ref{rank reuction } shows that $T$ is injective on the set of all rank-1 matrices. Let $E_{ij} \in \mathbb{E}$ and $C$ be a rank-1 matrix such that $T(C)=E_{ij}$. Let $C_{lk}\neq 0$. Then $C_{lk} E_{lk}\leq C$ implies that $C_{lk}T(E_{lk}) \leq E_{ij}$, which implies that $T(E_{lk})=\alpha E_{ij}$. If $C$ has more than one nonzero entry, say $C_{rs} \neq 0$ for $(l,k) \neq (r,s)$, then similarly we get that $T(E_{rs})=\beta E_{ij}$. Thus $T(\alpha E_{rs})=\alpha \beta E_{ij}=T(\beta E_{ij})$. But $\alpha  E_{rs} \neq  \beta E_{lk}$ leads to the contradiction that $T$ is injective on the set of all rank-1 matrices.

Thus $C$ is of the form $\alpha E_{lk}$. Now we show that $\alpha \in U(S)$. Since $T(C)=E_{ij}$ and $C=\alpha E_{lk}$, it follows that $\alpha T(E_{lk})=E_{ij}$. Thus there exists $\gamma \in S$ such that $\alpha \gamma =1$, which shows that $\alpha \in U(S)$. Recall that $C$ is the pre-image of $E_{ij}$. Since $T$ is injective, $E_{ij}$ can not have more than one pre-image. Since $\mathbb{E}$ is a finite set, we conclude that $T$ permutes $\mathbb{E}$ with some unit scalar multiplication. Hence $T$ is invertible. Part $(b)$ follows trivially.
\end{proof}
\end{corollary}

Suppose $T$ is an invertible linear operator on $M_{m \times n}(S)$. We know, by Corollary \ref{invertible of t}, that $T(E_{ij})=\alpha_{ij}E_{pq}$, for some $\alpha _{ij} \in U(S)$. We call $p$ and $q$ as first and second coordinates of $T(E_{ij})$ respectively and this p,q depend on i and j. For this invertible linear operator $T$, define the $m \times n$ array $\tau$ whose $(i,j)$th entry is $\tau(i,j)=\alpha_{ij}(p,q)$ for all $(i,j) \in \Delta$, where $\alpha _{ij}$ are units of $S$. The array $\tau$ is called the representation of $T$.

The following result can be proved similarly to the Lemma 3.7  \cite{Bapat} and so we omit the proof..
\begin{lemma}       \label{main lemma}
If $T$ is an invertible linear operator on $M_{m \times n}(S)$ that preserves the rank of every rank-1 matrix, then there exist permutations $\rho$ and $\sigma$ of $1,2, \cdots ,m$ and $1,2,\cdots ,n$ respectively such that
\begin{enumerate}
\item[(i)] $T(E_{ij})=\alpha _{ij}E_{\rho(i),\sigma(j)}$  or
\item[(ii)] $m=n$ and $T(E_{ij})=\alpha _{ij}E_{\sigma(j),\rho(i)}$ for all $(i,j) \in \Delta$, where $\alpha _{ij} \in U(S)$.
\end{enumerate}
\end{lemma}

\begin{lemma}       \label{u-v operator}
If $T$ satisfies the conclusions of Lemma \ref{main lemma}, then we have the following results.
\begin{enumerate}
\item[(i)] For all $i,l \in \{1,2, \cdots ,m\}$ and $j,k \in \{1,2, \cdots ,n\}$,
\begin{center}
$\alpha _{ij}\alpha _{lk}=\alpha _{lj}\alpha _{ik}$.
\end{center}

Thus there exist two diagonal matrices $C$ and $D$ such that $\alpha _{ij}=C_{ii}D_{jj}$ for all $i,l \in \{1,2, \cdots ,m\}$ and $j,k \in \{1,2, \cdots ,n\}$.
\item[(ii)] For any $ m \times n$ matrix $A$, there exist invertible matrices $U,V$ such that

$T(A)=UCADV$ if $m \neq n$ or $T(A)=VDA^tCU$ if $m=n$.
\end{enumerate}
\begin{proof}
(i) Since $E_{ij}+E_{ik}+E_{lj}+E_{lk}$ is a rank one matrix, $T(E_{ij}+E_{ik}+E_{lj}+E_{lk})$ is also a rank one matrix. By Lemma \ref{main lemma}, we have
\begin{center}
$T(E_{ij}+E_{ik}+E_{lj}+E_{lk})=\alpha_{ij}E_{\rho(i),\sigma(j)}+
\alpha_{ik}E_{\rho(i),\sigma(k)}+\alpha_{lj}E_{\rho(l),\sigma(j)}+\alpha_{lk}E_{\rho(l),\sigma(k)}$
\end{center}
Since $T$ preserves the rank of all rank one matrices, the rank of the matrix in the right hand side of the above equation is one. Hence
\begin{center}
$\alpha_{ij}\alpha_{lk}=\alpha_{lj}\alpha_{ik}$.
\end{center}

Take two diagonal matrices $C$ of order $m \times m$ and $D$ of order $n \times n$ defined by $C_{ii}=\alpha_{i1}\alpha ^{-1}_{11}$ and $D_{jj}=\alpha_{1j}$ for all $i \in \{1,2,3, \cdots ,m\}$ and $j \in \{1,2, \cdots ,n\}$. Then both $C$ and $D$ are invertible and $C_{ii}D_{jj}=\alpha _{ij}$.

(ii)Let $\pi$ be any permutation of $\{1,2, \cdots ,k\}$. Let $E^{m,n}_{ij}$ denote $m \times n$ matrix whose $(i,j)$th entry is $1$. Let $P_k(\pi)= \sum^k_{l=1}E^{k,k}_{l,\pi (l)}$. Then $P_k(\pi)$ is a permutation matrix of order $k$. But $E^{m,n}_{i,j}E^{n,r}_{u,v}=\delta _{j,u}E^{m,r}_{i,v}$, where $\delta_{j,u}$ is the Kroneker delta. Thus $E^{m,n}_{i,j}P_n(\pi)=E^{m,n}_{i,\pi(j)}$ and hence $P_m(\rho^{-1})E^{m,n}_{i,j}P_n(\sigma)=E^{m,n}_{\rho(i),\sigma(j)}$

If conclusion $(i)$ of Lemma \ref{main lemma} holds, then we define $U=P_m(\rho^{-1})$ and $V=P_n(\sigma)$. If $A$ is any $m\times n$ matrix, we have
\begin{align*}
T(A)=& T(\sum a_{ij}E^{m,n}_{ij})
\\=& \sum a_{ij}T(E^{m,n}_{ij})
\\=& \sum a_{ij} \alpha_{ij}E_{\rho(i),\sigma(j)}
\\=& \sum P_m(\rho^{-1})a_{ij} \alpha_{ij}E_{i,j}P_n(\sigma)
\\=& P_m(\rho^{-1}) \left(
                      \begin{array}{cccc}
                        a_{11}\alpha_{11} & a_{12}\alpha_{12} & \cdots & a_{1n}\alpha_{1n} \\
                        a_{21}\alpha_{21} & a_{22}\alpha_{22} & \cdots & a_{2n}\alpha_{2n} \\
                        \vdots & \vdots & \vdots & \vdots \\
                        a_{m1}\alpha_{m1} & a_{m2}\alpha_{m2} & \cdots & a_{m n}\alpha_{m n} \\
                      \end{array}
                    \right)
P_n(\sigma)
\\=& P_m(\rho^{-1})\left(
                      \begin{array}{cccc}
                        a_{11}C_{11}D_{11} & a_{12}C_{12}D_{12} & \cdots & a_{1n}C_{1n}D_{1n} \\
                        a_{21}C_{21}D_{21} & a_{22}C_{22}D_{22} & \cdots & a_{2n}C_{2n}D_{2n} \\
                        \vdots & \vdots & \vdots & \vdots \\
                        a_{m1}C_{m1}D_{m1} & a_{m2}C_{m2}D_{m2} & \cdots & a_{m n}C_{m n}D_{m n} \\
                      \end{array}
                    \right)
P_n(\sigma)
\\=&P_m(\rho^{-1})\left(
                      \begin{array}{cccc}
                        C_{11} & 0 & \cdots & 0 \\
                           0 & C_{22} & \cdots & 0 \\
                        \vdots & \vdots & \vdots & \vdots \\
                        0 & 0 & \cdots & C_{m m} \\
                      \end{array}
                    \right)\left(
                      \begin{array}{cccc}
                        a_{11} & a_{12} & \cdots & a_{1n} \\
                        a_{21} & a_{22} & \cdots & a_{2n} \\
                        \vdots & \vdots & \vdots & \vdots \\
                        a_{m1} & a_{m2} & \cdots & a_{m n} \\
                      \end{array}
                    \right)\left(
                      \begin{array}{cccc}
                        D_{11} & 0 & \cdots & 0 \\
                        0 & D_{22} & \cdots & 0 \\
                        \vdots & \vdots & \vdots & \vdots \\
                        0 & 0 & \cdots & D_{n n} \\
                      \end{array}
                    \right)
P_n(\sigma)
\\=& P_m(\rho^{-1})CAD P_n(\sigma)
\\=& UCADV
\end{align*}

Similarly, if conclusion $(ii)$ of Lemma \ref{main lemma} holds, then there exist $U=P^t_n(\rho^{-1})$ and $V=P^t_n(\sigma)$ such that
\begin{center}
$T(A)=VDA^tCU$
\end{center}

Thus $T$ is a $(U,V)$ operator.
\end{proof}
\end{lemma}
\begin{theorem}     \label{invertibility and u-v operator}
If $T$ is a linear operator on $M_{m \times n}(S)$, then the following statements are equivalent.
\begin{enumerate}
\item[(i)] $T$ is invertible and preserves the rank of all rank-1 matrices.
\item[(ii)] $T$ preserves the rank of all rank-1 and rank-2 matrices and preserves the dimension of all rank-1 subsemimodules.
\item[(iii)] $T$ is a $(U,V)$ operator.
\end{enumerate}
\begin{proof}
Lemma \ref{main lemma} and Lemma \ref{u-v operator} show that $(i)$ implies $(iii)$. Corollary \ref{rank-1 rank-2 invertible} shows that $(ii)$ implies $(i)$. To show that $(iii)$ implies $(ii)$, note that $(U,V)$ operators are always invertible, in fact $T^{-1}(A)=U^{-1}AV^{-1}$ or $T^{-1}(A)=U^{-1}A^tV^{-1}$. Such operators are preservers of all rank. The rest is implied by Lemma \ref{injectivity}.
\end{proof}
\end{theorem}
\section{Rank preservers of matrices over semirings}   \label{preservers of matrices}
In this section we characterize the linear operators those preserve the rank of matrices of any order over semirings.

We say that a linear operator $T$ on $M_{m \times n}(S)$ is a rank preserver if $T$ preserves the rank of all matrices.
\begin{lemma}      \label{decrease of rank}
Let $A,B \in M_{m \times n}(S)$ be two distinct matrices such that $r(A)=r(B)=1$, where $m>1,n>1$,
\begin{enumerate}
\item[(i)] If the number of nonzero entries in A is more than that of B, then there exists $C\in M_{m \times n}(S)$ such that $r(A+C)=1$ and $r(B+C)=2$.
\item[(ii)] If the number of nonzero entries in A is equal to that of B, then there exists $C\in M_{m \times n}(S)$ such that $r(A+C)=1$ and $r(B+C)=2$, or $r(A+C)=2$ and $r(B+C)=1$.
\end{enumerate}
\begin{proof}
(i) If $r(A+B)=2$, then it holds with $C=A$. So we assume that $r(A+B)=1$. Since the number of nonzero entries in A is more than that of B, there exists $(i_0,j_0) \in \Delta$ such that $A_{i_0j_0} \neq 0$ but $B_{i_0j_0}=0$. Consider the following cases:

\textbf{Case 1}: $A+B$ has at least two nonzero rows and and two nonzero columns.

Define $C$ by:
\begin{equation}
C_{ij}=
\left \{
\begin{array}{cc}
0 & \mbox{if $(i,j)=(i_0,j_0)$}\\
A_{ij}, & \mbox{otherwise}
\end{array}
\right.
\end{equation}

Then $A+C=A$ and hence $r(A+C)=1$. But $B+C$ is same as $B+A$ except $(i_0,j_0)$th entry. That is $(B+C)_{i_0j_0}=0$ but $(B+A)_{i_0j_0}=A_{i_0j_0} \neq 0$. We show that $r(B+C)=2$. Since $r(A+B)=1, A+B=a_{m \times 1}x^t_{1 \times n}$. If possible, let $B+C=b_{m \times 1}y^t_{1 \times n}$. Then $b_{i_0}y_{j_0}=0$ which implies that either $b_{i_0}=0$ or $y_{j_0}=0$. If $b_{i_0}=0$, all entries of $i_0$th row of $B+C$ are zero. So, B+C looks like $\left(
                                                   \begin{array}{c}
                                                     b_1 \\
                                                     \vdots \\
                                                     0 \\
                                                     \vdots \\
                                                     b_m \\
                                                   \end{array}
                                                 \right) \left(
                                                           \begin{array}{ccccc}
                                                             y_1 & \cdots & y_{j_0} & \cdots & y_n \\
                                                           \end{array}
                                                         \right)$

Since $B+C$ is same as $B+A$ except $(i_0,j_0)$th entry, $A+B$ looks like $\left(
                                                   \begin{array}{c}
                                                     a_1 \\
                                                     \vdots \\
                                                     a_{i_0} \\
                                                     \vdots \\
                                                     a_m \\
                                                   \end{array}
                                                 \right) \left(
                                                           \begin{array}{ccccc}
                                                             0 & \cdots & x_{j_0} & \cdots & 0 \\
                                                           \end{array}
                                                         \right)$

which shows that A+B has exactly one column that contradicts that A+B has at least two nonzero column.
Now B+C can be expressed as $\left(
                              \begin{array}{cc}
                                a_1 & a_1 \\
                                a_2 & a_2 \\
                                \vdots & \vdots \\
                                0 & a_{i_0} \\
                                \vdots & \vdots \\
                                a_m & a_m \\
                              \end{array}
                            \right)\left(
                                     \begin{array}{ccccc}
                                       x_1 & \cdots & x_{j_0} & \cdots & x_n \\
                                       x_1 & \cdots & 0 & \cdots & x_n \\
                                     \end{array}
                                   \right)$

So $r(B+C)=2$. If $y_{j_0}=0$, then proceeding in similar way, we can see that $A+B$ has exactly one nonzero row, which contradicts that $A+B$ has at least two nonzero row.

\textbf{Case 2}: $A+B$ has exactly one nonzero row, say $i_0$th row. Then the $i_0$th row of $A$ and $B$ are nonzero respectively, and all other rows are zero. Let $i_1$th row of $A$ be zero.

Define $C$ by:
\begin{equation}
C_{ij}=
\left \{
\begin{array}{cc}
0, & \mbox{if $(i,j)=(i_0,j_0)$}\\
A_{i_0j}+B_{i_0j}, & \mbox{if $i=i_0,j \neq j_0$}\\
A_{i_0j}+B_{i_0j}, & \mbox{if $i=i_1$}\\
0, & \mbox{otherwise}
\end{array}
\right.
\end{equation}

Then $A+C=\left(
            \begin{array}{cccccc}
              0 & \cdots & \cdots & \cdots & \cdots & 0 \\
              \vdots & \vdots & \vdots & \vdots &\vdots & \vdots \\
              A_{i_01}+B_{i_01} & A_{i_02}+B_{i_02} & \cdots & A_{i_0j_0} & \cdots & A_{i_0n}+B_{i_0n} \\
              \vdots & \vdots & \vdots & \vdots &\vdots & \vdots \\
              A_{i_01}+B_{i_01} & A_{i_02}+B_{i_02} & \cdots & A_{i_0j_0} & \cdots & A_{i_0n}+B_{i_0n} \\
              \vdots & \vdots & \vdots & \vdots &\vdots & \vdots \\
              0 & 0 & 0 & 0 & 0 & 0\\
            \end{array}
          \right)$
is of rank-1, but $B+C=\left(
            \begin{array}{cccccc}
              0 & \cdots & \cdots & \cdots & \cdots & 0 \\
              \vdots & \vdots & \vdots & \vdots &\vdots & \vdots \\
              A_{i_01}+B_{i_01} & A_{i_02}+B_{i_02} & \cdots & 0 & \cdots & A_{i_0n}+B_{i_0n} \\
              \vdots & \vdots & \vdots & \vdots &\vdots & \vdots \\
              A_{i_01}+B_{i_01} & A_{i_02}+B_{i_02} & \cdots & A_{i_0j_0} & \cdots & A_{i_0n}+B_{i_0n} \\
              \vdots & \vdots & \vdots & \vdots &\vdots & \vdots \\
              0 & 0 & 0 & 0 & 0 & 0\\
            \end{array}
          \right)$, which is same as $A+C$ except $(i_0,j_0)$th entry. That is $(B+C)_{i_0j_0}=0 $ but $(A+C)_{i_0j_0}=A_{i_0j_0} \neq 0$. This is similar to the case-1, so $r(B+C)=2$.

\textbf{Case 3}: $A+B$ has exactly one nonzero column. This is similar to Case 2.

(ii) If $r(A+B)=2$, then it holds with $C=A$. So we assume that $r(A+B)=1$. Since $A \neq B$, let us assume there exists $(i_0,j_0) \in \Delta$ such that $A_{i_0j_0} \neq B_{i_0j_0}$ and neither $A_{i_0j_0}$ nor $B_{i_0j_0}$ is zero. So either $A_{i_0j_0} \nleq B_{i_0j_0}$ or $B_{i_0j_0} \nleq A_{i_0j_0}$

\textbf{Case 1}: $A+B$ has at least two nonzero rows and and two nonzero columns.

\textbf{Sub-Case 1}: If $A_{i_0j_0} \nleq B_{i_0j_0}$, define $C$ by:

\begin{equation}
C_{ij}=
\left \{
\begin{array}{cc}
0, & \mbox{if $(i,j)=(i_0,j_0)$}\\
A_{ij}, & \mbox{otherwise}
\end{array}
\right.
\end{equation}
Then $A+C=A$, so $r(A+C)=1$ and $B+C$ is same as $A+B$ except $(i_0,j_0)$th entry. That is $(B+C)_{i_0j_0}=B_{i_0j_0}$ but $(A+B)_{i_0j_0}=A_{i_0j_0}+B_{i_0j_0} \neq B_{i_0j_0}$. We show that $r(B+C)=2$. Since $r(A+B)=1$, $A+B=a_{m \times 1}x^t_{1 \times n}$. If possible, let $B+C=b_{m \times 1}y^t_{1 \times n}$. So $a_{i_0}x_{j_0}=A_{i_0j_0}+B_{i_0j_0}=A_{i_0j_0}+b_{i_0}y_{j_0}$. Since $A+B$ has two nonzero rows and two nonzero columns, let other nonzero row and column be $i$th  row and $j$th column respectively. Then $a_i \neq 0$ and $x_j \neq 0$. Thus we get,
\begin{center}
$a_ix_j=b_iy_j$ \\ and
$a_{i_0}x_j=b_{i_0}y_j$
\end{center}
which implies that $a_ib_{i_0}=a_{i_0}b_i$, by multiplicative cancellative property of $S$.

Thus from $a_{i_0}x_{j_0}=A_{i_0j_0}+B_{i_0j_0}=A_{i_0j_0}+b_{i_0}y_{j_0}$ and $a_ix_{j_0}=b_iy_{j_0}$, we get
\begin{align*}
& a_ia_{i_0}x_{j_0}=a_i(A_{i_0j_0}+b_{i_0}y_{j_0})=a_{i_0}b_iy_{j_0}=a_ib_{i_0}y_{j_0}
\\ \Rightarrow & a_i(A_{i_0j_0}+b_{i_0}y_{j_0})=a_ib_{i_0}y_{j_0}
\\ \Rightarrow & A_{i_0j_0}+ b_{i_0}y_{j_0}=b_{i_0}y_{j_0}
\\ \textrm{i.e.}\;\; & A_{i_0j_0}+ B_{i_0j_0}=B_{i_0j_0}
\end{align*}
which contradicts that $A_{i_0j_0} \nleq B_{i_0j_0}$. So $r(B+C) \neq 1$.

Now, $B+C$ can be expressed as $\left(
                                 \begin{array}{cc}
                                   a_1 & a_1x_{j_0} \\
                                   a_2 & a_2x_{j_0} \\
                                   \vdots & \vdots \\
                                   a_{i_0} & B_{i_0j_0} \\
                                   \vdots & \vdots \\
                                   a_m & a_mx_{j_0} \\
                                 \end{array}
                               \right) \left(
                                         \begin{array}{cccccc}
                                           x_1 & x_2 & \cdots & 0 & \cdots & x_n \\
                                           0 & 0 & \cdots & 1 & \cdots & 0 \\
                                         \end{array}
                                       \right)$

So, $r(B+C)=2$.

\textbf{Sub-Case 2}: If $B_{i_0j_0} \nleq A_{i_0j_0}$, define $C$ by:

\begin{equation}
C_{ij}=
\left \{
\begin{array}{cc}
0, & \mbox{if $(i,j)=(i_0,j_0)$}\\
B_{ij}, & \mbox{otherwise}
\end{array}
\right.
\end{equation}
Then proceeding as above we get $r(A+C)=2$ and $r(B+C)=1$.

\textbf{Case 2}: $A+B$ has exactly one nonzero row, say $i_0$th row. Then the $i_0$th row of $A$ and $B$ are nonzero respectively, and all other rows are zero. Let $i_1$th row of $A$ be zero.

\textbf{Sub-Case 1}: If $A_{i_0j_0} \nleq B_{i_0j_0}$.
Define $C$ as the following:
\begin{equation}
C_{ij}=
\left \{
\begin{array}{cc}
B_{i_0j_0}, & \mbox{if $(i,j)=(i_0,j_0)$}\\
A_{i_0j}+B_{i_0j}, & \mbox{if $i=i_0,j \neq j_0,j_1$}\\
A_{i_0j}+B_{i_0j}, & \mbox{if $i=i_1,j \neq j_1$}\\
A_{i_0j_1}+B_{i_0j_1}+B_{i_0j_0}, & \mbox{if $i=i_0,i_1$ and $j=j_1$}\\
0, & \mbox{otherwise}
\end{array}
\right.
\end{equation}

Then $A+C=\left(
            \begin{array}{ccccccc}
              0 & \cdots & \cdots & \cdots & \cdots & \cdots & 0 \\
              \vdots & \vdots & \vdots & \vdots &\vdots & \vdots & \vdots \\
              A_{i_01}+B_{i_01} & \cdots & A_{i_0j_0}+B_{i_0j_0} & \cdots & A_{i_0j_1}+B_{i_0j_1}+B_{i_0j_0} & \cdots & A_{i_0n}+B_{i_0n} \\
              \vdots & \vdots & \vdots & \vdots &\vdots & \vdots & \vdots \\
              A_{i_01}+B_{i_01} &  \cdots & A_{i_0j_0}+B_{i_0j_0} & \cdots & A_{i_0j_1}+B_{i_0j_1}+B_{i_0j_0} & \cdots & A_{i_0n}+B_{i_0n} \\
              \vdots & \vdots & \vdots & \vdots &\vdots & \vdots & \vdots \\
              0 & 0 & 0 & 0 & 0 & 0 & 0\\
            \end{array}
          \right)$
is a matrix of rank-1, but

$B+C=\left(
            \begin{array}{ccccccc}
              0 & \cdots & \cdots & \cdots & \cdots & \cdots & 0 \\
              \vdots & \vdots & \vdots & \vdots &\vdots & \vdots & \vdots \\
              A_{i_01}+B_{i_01} & \cdots & B_{i_0j_0} & \cdots & A_{i_0j_1}+B_{i_0j_1}+B_{i_0j_0} & \cdots & A_{i_0n}+B_{i_0n} \\
              \vdots & \vdots & \vdots & \vdots &\vdots & \vdots & \vdots \\
              A_{i_01}+B_{i_01} &  \cdots & A_{i_0j_0}+B_{i_0j_0} & \cdots & A_{i_0j_1}+B_{i_0j_1}+B_{i_0j_0} & \cdots & A_{i_0n}+B_{i_0n} \\
              \vdots & \vdots & \vdots & \vdots &\vdots & \vdots & \vdots \\
              0 & 0 & 0 & 0 & 0 & 0 & 0\\
            \end{array}
          \right)$, which is same as $A+C$ except $(i_0,j_0)$th entry. That is $(B+C)_{i_0j_0}=B_{i_0j_0} $ but $(A+C)_{i_0j_0}=A_{i_0j_0}+B_{i_0j_0} \neq B_{i_0j_0}$. This is similar to the Sub-Case 1 of Case 1, so $r(B+C)=2$.

\textbf{Sub-Case 2}: If $B_{i_0j_0} \nleq A_{i_0j_0}$. Define $C$ by:
\begin{equation}
C_{ij}=
\left \{
\begin{array}{cc}
A_{i_0j_0}, & \mbox{if $(i,j)=(i_0,j_0)$}\\
A_{i_0j}+B_{i_0j}, & \mbox{if $i=i_0,j \neq j_0,j_1$}\\
A_{i_0j}+B_{i_0j}, & \mbox{if $i=i_1,j \neq j_1$}\\
A_{i_0j_1}+B_{i_0j_1}+B_{i_0j_0}, & \mbox{if $i=i_0,i_1$ and $j=j_1$}\\
0, & \mbox{otherwise}
\end{array}
\right.
\end{equation}
Then proceeding as above we get $r(A+C)=2$ and $r(B+C)=1$.

\textbf{Case 3}: $A+B$ has exactly one nonzero column. This is similar to Case 2.
\end{proof}
\end{lemma}
\begin{lemma}      \label{decrease of rank of matrix}
If $T$ is a linear operator on $M_{m \times n}(S)$ with $m > 1, n > 1$ and $T$ is not invertible but preserves the rank of all rank-1 matrices. Then $T$ decreases the rank of some rank-2 matrix to 1.
\begin{proof}
By the proof of Corollary \ref{rank-1 rank-2 invertible}, $T$ is not injective on the set of all rank-1 matrices. So there exist distinct rank-1 matrices X and Y such that $T(X)=T(Y)$. Without lose of generality we may assume that the number of nonzero entries in $X$ is more than or equal to that of $Y$. Then from Lemma \ref{decrease of rank}, there exists a matrix $C$ such that
\begin{align*}
&\textrm{either}\;\; r(X+C)=2, r(Y+C)=1
\\&\textrm{or} \;\; r(X+C)=1, r(Y+C)=2
\end{align*}
Thus $T(X+C)=T(X)+T(C)=T(Y)+T(C)=T(Y+C)$ is a rank-1 matrix. Thus $T$ decreases the rank of rank-2 matrix $X+C$. Similarly, in the latter case, $T$ decreases the rank of $Y+C$.
\end{proof}
\end{lemma}
\begin{theorem}  \label{rank preserver}
Let $T$ be a linear operator on $M_{m \times n}(S)$ where $m > 1, n > 1$. Then $T$ is a rank preserver if and only if $T$ is a $(U,V)$ operator.
\begin{proof}
From Theorem \ref{invertibility and u-v operator} and Lemma \ref{decrease of rank of matrix}, we see that the necessity of the condition is satisfied. The sufficiency is trivial, since every $(U,V)$ operator is a rank preserver.
\end{proof}
\end{theorem}
\begin{theorem}
Let $T$ be a linear operator on $M_{m \times n}(S)$ where $m>1,n>1$. Then $T$ is a rank preserver if and only if $T$ preserves the rank of all rank-1 and rank-2 matrices.
\begin{proof}
If $T$ preserves the rank of all rank-1 and rank-2 matrices, then $T$ is invertible by Lemma \ref{decrease of rank of matrix}. Thus $T$ is a rank preserver by Theorem \ref{invertibility and u-v operator} and Theorem \ref{rank preserver}.

The converse is trivial.
\end{proof}
\end{theorem}

\bibliographystyle{amsplain}

\end{document}